\documentclass{amsart}
\usepackage[normalem]{ulem}
\numberwithin{equation}{section}
\usepackage{longtable}
\usepackage[all]{xy}
\usepackage{amssymb}
\usepackage{comment}
\usepackage{stmaryrd}
\usepackage{xcolor}
\usepackage{mathtools}
\usepackage{enumitem}
\usepackage[OT2]{fontenc}
\usepackage[T1]{fontenc}
\definecolor{Vino}{rgb}{0.356,0,0}
\definecolor{Ruta}{rgb}{0.309, 0.459, 0.208}
\usepackage[final,colorlinks,linkcolor=gray,citecolor=Vino,urlcolor=black]{hyperref}
\let\cal\mathcal
\def\Ascr{{\cal A}}
\def\Mscr{{\cal B}}
\def\Bscr{{\cal B}}
\def\Cscr{{\cal C}}
\def\Dscr{{\cal D}}
\def\Escr{{\cal S}}
\def\Eescr{{\cal E}}

\def\Hscr{{\cal H}}

\def\Kscr{{\cal K}}

\def\Mscr{{\cal M}}

\def\Oscr{{\cal O}}
\def\Pscr{{\cal P}}

\def\Tscr{{\cal T}}

\let\blb\mathbb
\def\CC{{\blb C}}

\def \PP{{\blb P}}

\def \ZZ{{\blb Z}}

\def \RR{{\blb R}}

\def\id{\text{id}}

\def\quot{/\!\!/}

\def\mod{\operatorname{mod}}

\def\gr{\operatorname{gr}}

\def\coh{\mathop{\text{\upshape{coh}}}}

\def\gr{\operatorname {gr}}
\def\Spec{\operatorname {Spec}}

\def\Hom{\operatorname {Hom}}
\def\uHom{\operatorname {\mathcal{H}\mathit{om}}}
\def\uExt{\operatorname {\mathcal{E}\mathit{xt}}}
\def\uEnd{\operatorname {\mathcal{E}\mathit{nd}}}

\def\id{{\operatorname {id}}}

\DeclareMathOperator{\Proj}{Proj}

\newtheorem{lemma}{Lemma}[section]
\newtheorem{proposition}[lemma]{Proposition}
\newtheorem{theorem}[lemma]{Theorem}
\newtheorem{corollary}[lemma]{Corollary}

\theoremstyle{definition}

\newtheorem{example}[lemma]{Example}
\newtheorem{definition}[lemma]{Definition}
\newtheorem{conjecture}[lemma]{Conjecture}

{

}

\theoremstyle{remark}

\newtheorem{remark}[lemma]{Remark}

\newdimen\uboxsep \uboxsep=1ex
\def\uboxn#1{\vtop to 0pt{\hrule height 0pt depth 0pt\vskip\uboxsep
\hbox to 0pt{\hss #1\hss}\vss}}

\def\uboxs#1{\vbox to 0pt{\vss\hbox to 0pt{\hss #1\hss}
\vskip\uboxsep\hrule height 0pt depth 0pt}}

\let\oldmarginpar\marginpar
\def\marginpar#1{\oldmarginpar{\raggedright \tiny \baselineskip 0pt \lineskip 0pt #1}}

\title[HMS symmetries and hypergeometric systems]{HMS symmetries and hypergeometric systems}
\author[\v{S}pela \v{S}penko]{\v{S}pela \v{S}penko} 
\address{\newline D\'epartement de Math\'ematique, Universit\'e Libre de Bruxelles, Campus de la Plaine CP 213, Bld du Triomphe, B-1050 Bruxelles\newline {\tt spela.spenko@ulb.be}}
\thanks{The work on this survey was partially supported by the FNRS project MIS/BEJ F.4545.21 and the ARC project ``Algebra''.}

\begin{document}
\begin{abstract}
The derived category of an algebraic variety might be a source of a myriad of new (categorical) symmetries. Some are predicted by homological mirror symmetry, to be obtained from the fundamental group of the space of complex structures of its mirror partner. These finally lead to differential equations. We expositorily unravel a part of this conjectural master plan for a class of toric varieties.
\end{abstract}

\maketitle

\hfill {\it nasvidenje, Marjan}

\hfill {\it\small nekoč... nekje...}

\section{Overview}
Hilbert's 21st problem asks about the existence of Fuchsian linear differential equations on the Riemann sphere with prescribed singular points and  monodromy representation of the fundamental group of the complement of the singular points \cite{Hilbert1900}. The first (slightly erroneous) solution was proposed by the Slovenian mathematician Plemelj \cite{Plemelj}. A suitably adapted version of this problem was solved and generalised, depending on the context, by Deligne \cite{Deligne},  Kashiwara \cite{Kashiwara}, Mebkhout \cite{Mebkhout}, Beilinson - Bernstein \cite{BB}, ... The solution is now known as the Riemann-Hilbert correspondence.

Homological mirror symmetry predicts the existence of an action of the fundamental group of the "stringy Kähler moduli space (SKMS)"  on the derived category of an algebraic variety. The prediction was established by Halpern-Leistner and Sam for certain toric varieties \cite{HLSam}. A decategorification of this action yields a representation of the fundamental group of the SKMS,   and our joint work with Michel Van den Bergh shows that it corresponds under the Riemann-Hilbert correspondence to a hypergeometric system of differential equations \cite{SVdBgkz}.

In this expository note we aim to explain the above terms and finally present the mentioned results. 

{
  \hypersetup{linkcolor=black}
  \tableofcontents
}

\section{Hilbert's 21st problem}
We begin with a classical problem, namely Hilbert's 21st problem. It is a part of 
the list of 23 problems \cite{Hilbert1900,Hilbert1902}, published by Hilbert in 1900, which has been  influential for the future mathematical development. The 21st one had the following formulation:

\medskip

{\it{To show that  there  always exists a linear  differential  equation of  the  Fuchsian class, with given singular points  and monodromic  group.}}

\medskip

We shall first decipher the problem a little bit. 

\subsection{Fuchsian type}
A system of linear differential equations
\begin{equation}\label{ode}
\left(
\begin{matrix}
y_1'\\
\vdots\\
y_n'
\end{matrix}
\right)
=A(z) 
\left(
\begin{matrix}
y_1\\
\vdots\\
y_n
\end{matrix}
\right)
\end{equation}
is of {\em Fuchsian} type if {$A(z)$ is holomorphic on $\bar{\CC}\setminus\{a_1,\dots,a_N\}$ with a pole of order $1$ at $a_j$, $1\leq j\leq N$, where we denote $\bar{\CC}=\CC\cup \{\infty\}$.

In particular,\footnote{We may take $y_i=y^{(i)}$, $0\leq i\leq n$ and for $A$ an ($(n+1)\times (n+1)$-)matrix with nonzero entries only  on the first upper diagonal where they are equal to $1$ and in the last row.} $\sum_{i=0}^n q_{i}(z)y^{(n-i)}=0$, $q_n(z)=1$, is Fuchsian if and only if the familiar Fuchsian condition is satisfied, i.e. $q_i(z)(z-a)^i$ is holomorphic at $z=a$ for $a\in \CC$ and $q_i(z)z^i$ is holomorphic at $z=\infty$, for  $0\leq i\leq n$.\footnote{This follows by taking the $n\times n$-matrix  with $a_{i,i+1}=-1$, $a_{n,i}=q_{n-i+1}$, and $a_{ij}=0$ otherwise.}

\subsection{Monodromy}
Assume that we have a system of linear differential equations \eqref{ode} with singularities at finitely many points $\{a_1,...,a_N\}$. Let $\gamma$ be a closed path (so $\gamma(0)=\gamma(1)$) in $\bar{\CC}\setminus \{a_1,...,a_N\}$.

Let $y_1,\dots, y_n$ be a basis of solutions of the system on an open set around $\gamma(0)$ (they exist by the local existence theorem for differential equations). These solutions are guaranteed to exist a priori only locally. However, we can analytically continue them along $\gamma$. Let us denote by $\tilde{y}_1,\dots, \tilde{y}_n$ analytic continuations of $y_1,\dots,y_n$ along $\gamma$. 

Because both $y_1,\dots,y_n$ and $\tilde{y}_1,\dots, \tilde{y}_n$ form a basis of solutions around $\gamma(0)$ they should be related via an invertible linear map. We denote it  by $\rho_{\gamma}$. It turns out that $\rho_{\gamma}$ only depends on the homotopy class of $\gamma$. 
Therefore we obtain a group homomorphism  
\[
\pi_1(\bar{\CC}\setminus  \{a_1,...,a_N\})\to {\rm GL}_n(\CC),\quad [\gamma]\mapsto\rho_\gamma.\] 
This is what is called a {\em monodromy} representation. 

\medskip
To close this discussion, we look at a concrete example of a differential equation.
\begin{example}
We take the differential equation 
$zy'-\alpha y=0$.  First note that it has singularities at $0$ and at $\infty$. (It is of Fuchsian type.)   
We take a loop $\gamma$ around $0$. A local solution is equal to $y=z^\alpha$, its analytic continuation along $\gamma$ equals $\tilde{y}=e^{2\pi i \alpha}z^{\alpha}$. To construct a monodromy representation we first notice that the fundamental group of $\bar{\CC} \setminus \{0, \infty\}=\CC\setminus\{0\}$,  is isomorphic to $\ZZ$, and we can identify the generator $1$ with the homotopy class of $\gamma$. The monodromy representation is then given by 
\[
\rho:\pi_1(\bar{\CC}\setminus \{0,1\})\cong \ZZ\to {\rm GL}_1(\CC)\cong \CC^*, \quad k\mapsto e^{2\pi \mathrm{i} k\alpha}.
\]
\end{example}

\subsection{Formulation}
Let us now restate the problem. As input we have
\begin{itemize}
	\item
a finite set of points $\{a_1,...,a_N\}$, and 
\item
a representation $\rho$ of $\pi_1(\bar{\CC}\setminus \{a_1,...,a_N\})$. 
% i.e.  a group homomorhism $\rho: \pi_1(\bar{\CC}\setminus \{a_1,...,a_N\})\to {\rm GL}_n(\CC)$ for some $n\in \NN$. 
\end{itemize}
Then Hilbert's 21st problem reads as follows:
{\it Does there exist a system of linear differential equations of Fuchsian type with singular points $\{a_1,\dots,a_N\}$ and the monodromy representation equal to $\rho$?}

\subsection{Progress}

Already in 1908  Plemelj proposed a complete solution \cite{Plemelj}. %\footnote{With the congress have taken place in Slovenia let me add some historic facts. Plemelj was one of the first Slovenian mathematicians, and still around one fourth of Slovenian mathematicians are related to him (by math genealogy).%He was also the first rector of the University of  Ljubljana which will in 2021 turn 102 years old.
%}
Unfortunately, it turned out that Plemelj's solution was not entirely correct. (Nevertheless, Plemelj's proof shows that one can find a system of linear differential equations which is Fuchsian at all but one point, where it is regular, see \S\ref{sec:reg} below.)  
In 1988 Bolibrukh found a counterexample for $N=4$ and a $\rho$ of degree $3$ \cite{Bolibrukh}.

 The problem  then transformed into classifying the input data that correspond to systems of differential equations of Fuchsian type.

\medskip
Among algebraic geometers the focus was however directed towards  higher dimensions with a suitably rendered condition. Instead of Fuchsian type one requires regularity, a weaker condition.

\section{Riemann-Hilbert correspondence} 
%We first present a formulation and solution of the Hilbert's 21st problem in higher dimensions, and then its generalisation in the context of D-modules. 
There are plentiful variants of the Riemann-Hilbert correspondence. We first present one in line with the previous discussion, and then its powerful generalisation to the context of D-modules. We mostly follow \cite{hottabook}. We also mention \cite{Katz} for a very nice review of Deligne's work on  Hilbert's 21st problem.

\subsection{Integrable connections}
We first need to make sense of differential equations on general manifolds where we have no global coordinates at our disposal. 

Let $X$ be a complex manifold. Let $\Tscr_X$ be the tangent sheaf on $X$ (i.e. the sheaf of vector fields).\footnote{Note that $\Tscr_X$ may also be identified with derivations in $\uEnd_{\CC_X}(\Oscr_X)$.} 
\begin{definition}
An integrable connection on $X$ is a pair $(M,\nabla)$ where $M$ is finite dimensional vector bundle on $X$ and a linear map $\nabla:\Tscr_X\otimes \tilde{M}\to \tilde{M}$ where $\tilde{M}$ is the sheaf of sections of $M$ such that\footnote{We use  standard notation $\nabla_\theta(m):=\nabla(\theta\otimes m)$.}
\begin{itemize}
	\item $\nabla_{f\theta}(m)=f\nabla_\theta(m)$ for $f\in \Oscr_X$, $\theta\in \Tscr_X$, $m\in \tilde{M}$,
	\item $\nabla_\theta(fm)=\theta(f)m+f\nabla_\theta(m)$ for $f\in \Oscr_X$, $\theta\in \Tscr_X$, $m\in \tilde{M}$,
	\item  $\nabla_{[\theta_1,\theta_2]}(m)=[\nabla_{\theta_1},\nabla_{\theta_2}](m)$ for $\theta_1,\theta_2\in \Tscr_X$, $m\in \tilde{M}$.
\end{itemize}
\end{definition}

With the natural definition of morphisms we obtain an abelian category of connections on $X$ which we denote by  ${\rm Conn}(X)$.

\begin{remark}\label{rem:decon}
For a system of differential equations \eqref{ode} on  $X=\CC\setminus \{0\}$ (i.e. $0$ is the only singularity different from $\infty$), $M$ is the trivial vector bundle of rank $n$, and $\nabla$ is given by $\nabla_{\partial/\partial z}(y)=y'-A(z)y$ for $y\in \tilde{M}=(\Oscr_X)^n$.

Conversely, if $(M,\nabla)$ is an integrable connection on $X=\CC\setminus \{0\}$ then $M$ is a trivial vector bundle, say of rank $n$. We choose an $\Oscr_X$-basis $(e_i)_i$ of $\tilde{M}= \Oscr_X^n$. Define $a_{ij}(z)$, $1\leq i,j\leq n$, by $\nabla_{\partial/\partial z}(e_j)=-\sum_{i=1}^n a_{ij}(z)e_i$. Then $\nabla_{\partial/\partial z}(y)=\nabla_{\partial/\partial z}(\sum_i y_i e_i)=\sum_i y_i'e_i+\sum_i y_i\nabla_{\partial/\partial z}(e_i)=y'-A(z)y$ for $y\in \tilde{M}$.
\end{remark}

The solutions of an integrable connection are defined as $\{m\in \tilde{M}\mid \nabla_\theta(m)=0\;\text{for all $\theta\in \Tscr_X$}\}$ and are called {\em horizontal sections}.

%\medskip

\subsection{Meromorphic connections}
We now extend the concept of integrable connections to allow poles as well. Let $D\subset X$ a divisor. 
Let $\Oscr_X[D]$ be a sheaf of meromorphic functions on $X$, holomorphic on $X\setminus D$ with poles along $D$.

\begin{definition}
A coherent $\Oscr_X[D]$-module $M$\footnote{We note  that the definition implies that the restriction $M_{X\setminus D}$ of a meromorphic connection $M$  to $X\setminus D$ is a locally free $\Oscr_{X\setminus D}$-module.} is a {\em meromorphic connection} if there exists a  map $\nabla:M\to \Omega^1_X\otimes_{\Oscr_X} M$ such that 
\begin{itemize}
	\item
$\nabla(fs)=df \otimes s +f \nabla s$,
\item
$[\nabla_{\theta},\nabla_{\theta'}]=\nabla_{[\theta,\theta']}$ for $\theta,\theta'\in \Tscr_X$ (where $\nabla_\theta:M\to M$ is $\nabla'_\theta$ for $\nabla': \Tscr_X\otimes M\to M$ obtained from $\nabla$).
\end{itemize}
\end{definition}

With the natural definition of morphisms between meromorphic connections, we obtain an abelian category  ${\rm Conn}(X;D)$ of meromorphic connections. 

\begin{remark}
This remark is an analogue of Remark \ref{rem:decon}. We obtain a natural one to one correspondence between linear differential equations on $\CC$ with possible poles at $0$ and meromorphic connections in ${\rm Conn}(X;D)$. 
\end{remark}

\subsection{Regular singularities}\label{sec:reg}
Here we define regular singularities of differential equations, which are a generalisation of the differential equations of Fuchsian type.

\begin{definition}\label{def:regular_classical}
In complex dimension $1$ a system of differential equations has {\em regular singularities} if every solution $y$ on a punctured angular sector around a singular point in $\{a_1,\dots,a_N\}$ has moderate growth, i.e.
\begin{itemize}
	\item
$a_j$ finite: $|y(z)|= O(|z-a_j|^{-m})$ for some $m\geq 0$ as $z\to a_j$,
\item
$a_j=\infty$: $|y(z)|= O(|z|^m)$ for some $m\geq 0$ as $z\to \infty$.
\end{itemize}
\end{definition}

This also has an algebraic interpretation which can be moreover generalised to higher dimensions and all manifolds. %Here we for simplicity give a definition that applies only in the case when singularities are contained in  a normal crossing divisor (i.e.  a divisor that  can be locally defined  by a function $x_1\cdots x_r$ where $(x_1,\dots,x_n)$ are local coordinates), but Proposition \ref{prop:curve_testing} below gives a characterisation that holds for general divisors. 

%Let $D$ be a normal crossing divisor.  Then $D$ is a union of $r$ smooth divisors $D_i$ which intersect transversally. Let $\Tscr_X[D]$ be the subsheaf of the tangent sheaf which preserves the ideals of each of the $D_i$. At a point $x$ where all $D_i$ intersect $\Tscr_X[D]$ is a free $\Oscr_X$-module on $x_i\partial/\partial x_i $, $1\leq i\leq r$, $\partial/\partial x_i$, $r< i\leq n$. 

\begin{definition}\label{def:regular}
A meromorphic connection $(M,\nabla)$  in ${\rm Conn}(X;D)$ is {\em regular} if $(i^*M)_0$ is regular for every $i:B\to X$ such that $i^{-1}D=\{0\}$. 
\end{definition}
%\begin{definition}
%Let $U$ be an open set in $X$ such that $D=X\setminus U$ is a normal crossing divisor. Let $(M,\nabla)$ be an integrable connection on $U$. Then $(M,\nabla)$ is {\em regular} on $(X,D)$ if there exist a vector bundle $N$ on $X$ which extends $M$ and ${\nabla}^D:\Tscr_X[D]\otimes\tilde{N}\to \tilde{N}$ which extends $\nabla$. 
%\end{definition}

\begin{comment}

It turns out that Definition \ref{def:regular} coincides with the classical Definition \ref{def:regular_classical} (cf. Remark \ref{rem:decon}). %Let $\tilde{K}$ denote the set of (possibly) multilvalued holomorphic functions defined on $B(0,\epsilon)\setminus \{0\}$.

\begin{proposition}\label{prop:reg}
An integrable connection $(M,\nabla)$ on $(\CC,\{0\})$ is regular at $x=0$ if and only if all horizontal sections of $\nabla$ have moderate growth.\footnote{As mentioned above this characterisation holds also for  general divisors, for definition of regular conection in this case we refer to \cite{Hotta}.}
\end{proposition}

We denote by $B$ the unit disk in $\CC$. The following proposition says that we may test regularity (locally) in dimension $1$. %(It holds for all divisor even though Definition \ref{def:regular} only applies to normal crossing divisors.)

\begin{proposition}\label{prop:curve_testing}
Let $D$ be a divisor in $X$. An integrable connection $(M,\nabla)$ on $U=X\setminus D$ is regular on $(X,D)$ iff $i^*M$ is regular for every $i:B\to X$ such that $i^{-1}D=\{0\}$. 
\end{proposition}
\end{comment}

We also mention that with the natural definition of morphisms between regular meromorphic connections on $(X,D)$ we obtain an abelian category  ${\rm Conn}^{reg}(X;D)$.

\subsection{Deligne's Riemann-Hilbert correspondence}\label{subsec:DRH}

\begin{theorem}\cite{Deligne}\label{thm:Deligne}
Let $X$ be a complex manifold and let $D$ be a divisor in $X$. Then the restriction functor induces an equivalence ${\rm Conn}^{reg}(X;D)\xrightarrow{\sim} {\rm Conn}(X\setminus D)$.%\footnote{Here ${\rm Conn}(X\setminus D)={\rm Conn}(X\setminus D,\emptyset)$.}.
\end{theorem}

Deligne's theorem constitutes the essential part of the correspondence between systems of differential equations on $X$ with regular singularities along $D$ and representations of the fundamental group of $X\setminus D$.  

\begin{corollary}\label{cor:Deligne}
There is an equivalence of categories between ${\rm Conn}^{reg}(X;D)$ and ${\rm rep}(\pi_1(X\setminus D))$. 
\end{corollary}

This equivalence factors as 
\begin{equation}\label{diagram}
%\left(
\xymatrix{
%{\rm Conn}^{reg}(X;D)\xrightarrow{\sim} {\rm Conn}(X\setminus D)\xrightarrow{\sim} {\rm Loc}(X\setminus D)\xrightarrow{\sim} \pi_1(X\setminus D)
{\rm Conn}^{reg}(X;D)\ar[r]^-{\sim}\ar[d]_-{\wr}&{\rm rep}(\pi_1(X\setminus D))\\
{\rm Conn}(X\setminus D)\ar[r]^-{\sim}&{\rm Loc}(X\setminus D)\ar[u]_-{\wr}}
%\right)
\end{equation}
 where ${\rm Loc}(X\setminus D)$ is the category of local systems, i.e. locally constant sheaves of finite dimensional $\CC$-vector spaces. 
The first (vertical) equivalence is the restriction from Theorem \ref{thm:Deligne}, the second is obtained by taking the horizontal sections (``solutions of the system'') and the last (vertical) arrow is a well-known equivalence (see e.g. \cite{Achar}) which sends a local system $L$ to the representation of $\pi_1(X\setminus D)$ on $L_{x_0}$ that associates to   every path  an isomorphism of $L_{x_0}$ along itself (which exists as $L$ is locally constant).

\medskip
The statement holds also in the context of smooth algebraic varieties which was Deligne's original motivation. 

In short, we could say that topology, here measured by the fundamental group, is somewhat determined by  analysis or algebra, here represented by differential equations with regular singularities. 

\subsection{D-modules}\label{sec:Dmodules}
We continue towards a  generalisation of Deligne's correspondence to other systems of linear differential equations. 

For this we move on the left hand-side of the above diagram a bit more towards the algebra direction, and replace the differential equations with modules over the ring of differential operators. We enter the framework of so-called D-modules. 
%The thread of thinking that we present here was instilled in us by reading 
We follow \cite[Introduction]{hottabook}.

\medskip

Let $X$ be an open submanifold in $\CC^n$ and let $\Oscr(X)$ be holomorphic functions globally defined on $X$. With $D$ we denote the set of partial differential operators with coefficients in $\Oscr(X)$. Namely, 
\[
D=\left\{\sum_{i_1,\dots,i_n} f_{i_1\dots i_n}\left(\frac{\partial}{\partial x_1}\right)^{i_1}\cdots \left(\frac{\partial}{\partial x_n}\right)^{i_n}\mid f_{i_1\dots i_n}\in \Oscr(X)\right\},
\]
where $x_i$ are coordinate functions on $\CC^n$. Note that $D$ also has a ring structure. For example, $D$ contains the $n$-th Weyl algebra for $X=\CC^n$ (we take only polynomial coefficients).

Now take $P$ in $D$. Then $P$ corresponds to a differential equation.\footnote{For example, $x\frac{\partial}{\partial x}-\alpha$ corresponds to the equation $xy'-\alpha y$=0.} We can represent the holomorphic (global) solutions as follows:
\[
\{u\in \Oscr(X)\mid Pu=0\}\cong \Hom_{D}(D/DP,\Oscr(X)),\quad u\mapsto (d\mapsto du).
\]
%One only needs to notice that any  homomorphism on the right-hand side is determined by the image $u$ of $1$ which for it to be well defined should satisfy $Pu=0$.

We can proceed similarly if we have a collection of $P_{ij}\in D$, $1\leq i \leq k$, $1\leq j\leq l$, corresponding to a system of differential equations. Then the solution $(u_j)$ of the system given by the matrix $(P_{ij})$ can be identified with 
\[
\{(u_j)\mid (P_{ij})(u_j)=0\}\cong \Hom_{D}(M,\Oscr(X))
\]
where $M$ is defined by the short exact sequence 
\[
D^k\xrightarrow{(P_{ij})} D^l\to M\to 0.
\] 

In sum, we have found a way to turn systems of differential equations into finitely presented $D$-modules, and have described their (global) solutions purely algebraically using homomorphisms.  

However, solutions may not exist globally, so therefore we should use a tool that takes into account also local solutions. From modules we should pass to sheaves, as we have already done in the beginning of this section. Now $\Oscr$ denotes the sheaf of holomorphic functions. 
% (for an open set $U\subset X$, $\Oscr(U)$ is the set of holomorphic functions on $U$). 
Similarly, we replace $D$ by $\Dscr$ ($\Dscr(U)$ consists of partial differential operators with coefficients in $\Oscr(U)$). Then we can look at the sheaf $\uHom_{\Dscr}(\Mscr, \Oscr)$ ($U\mapsto\Hom_{\Dscr(U)}(\Mscr(U),\Oscr(U))$). %From now on when we say a $D$-module we mean a $\Dscr$-module, so a sheaf of $D$-modules.

There is another caveat to consider. We may be interested in relating different systems of differential equations, i.e. from solutions of two systems deduce something about solutions of the system that is formed as the union of the two systems. The problem that we encounter here is that the functor $\uHom_\Dscr(-, \Oscr)$ is not exact. So we should also consider ``higher solutions'', namely the extension modules $\uExt^i_\Dscr(\Mscr,\Oscr)$. 

%\medskip

It will turn out that higher solutions give us almost all the topological data that we need. Perhaps it is then a good point to ask what kind of sheaves these higher solutions are.  We know they are sheaves of $\CC$-vector spaces. Is there any other property that distinguishes them?

Recall from \eqref{diagram} (applied with $D=\emptyset$) that if $\Mscr$ is associated to a connection  then we obtain a local system, i.e. a locally constant sheaf of finite dimensional $\CC$-vector spaces. It turns out that this correspondence generalises if we restrict 
to holonomic modules\footnote{A coherent $\Dscr_X$-module $M$ is holonomic if $\dim {\rm Ch}(M)=\dim X$. Here ${\rm Ch}(M)$ denotes the characteristic variety of $M$, i.e. the support of the associated graded module ${\gr M}$ (for a ``good'' filtration) on the cotangent bundle of $X$.}, they are those that roughly speaking give finite dimensional (higher) solution spaces. With this assumption, all the higher solution sheaves $\uExt^i_\Dscr(\Mscr,\Oscr)$ are {\em constructible}, which means that they are built from local systems. More precisely, there exists a stratification of $X=\sqcup_\alpha X_\alpha$ into locally closed sets such that $F_i|_{X_\alpha}$ is a local system for all $i$. 

This is a prelude to a  correspondence between holonomic $\Dscr$-modules on the algebraic side and constructible sheaves on the topological side. Note that on the topological side we obtain an entire sequence of constructible sheaves, and to compute those we should also know something about the projective resolution of the modules, again on the algebraic side. A convenient machinery to process all this data at once and without losing too much information is the derived category. 

\subsection{Derived categories}
Let $\Ascr$ be an abelian category, for example the category $\mod({\Dscr_X})$ of $\Dscr_X$-modules on $X$, or the category $\mod(\CC_X)$ of sheaves of finite dimensional vector spaces on $X$, the categories that we have just seen. 

Let $C(\Ascr)$ be the category of complexes on $\Ascr$. We say that a map $f:X^\bullet\to Y^\bullet$ between two complexes is a quasi-isomorphism if it induces isomorphisms on cohomology, i.e. $H^i(f):H^i(X^\bullet)\xrightarrow{\sim} H^i(Y^\bullet)$ for all $i$. 

We want that the derived category does not distinguish between two complexes which are connected via a quasi-isomorphism. So we formally invert quasi-isomorphisms (see e.g. \cite[04VB]{stacks-project} for localisation in categories) and define the derived category as 
\[
D(\Ascr)=C(\Ascr)[{\rm qis}^{-1}].
\]
% Moreover, the derived category also contains the intial abelian category identified with the complexes leaving in degree 0.

Furthermore, if a covariant, resp. contravariant, functor $F:\Ascr\to \Bscr$ between two abelian categories (with $\Ascr$ having  enough injectives, resp. projectives) is left-exact then there exists a corresponding functor $RF:D^+(\Ascr)\to D^+(\Bscr)$, resp. $RF:D^-(\Ascr)\to D^+(\Bscr)$, between the derived categories (of bounded below, resp. above - below,  complexes). 

Let us zoom this in on our example. 
\begin{example}
We take for $F$ the solution functor, $F=\uHom_\Dscr(-,\Oscr)$. Then the derived functor $RF:D^-(\Dscr_X)^{o}\to D^+(\CC_X)$ is such that its cohomology sheaves are exactly the higher solutions, i.e. $H^i(RF)=\uExt^i(-,\Oscr_X)$. So the derived solution functor carries the information about all higher solutions. (Note that here and later we for brevity omit writing $\mod$.) 
\end{example}

\subsection{Riemann-Hilbert correspondence}\label{subsec:RH}
We are ready to state the Riemann-Hilbert correspondence in its full power and complexity, to connect all the module data with the data of solutions and higher solutions. 

We need to restrict to a subclass of complexes of $\Dscr_X$-modules that have regular\footnote{For the definition of regularity for $\Dscr_X$-modules on a complex manifold $X$ we refer to \cite[Definition 6.1.8]{hottabook}. %We note that a curve testing criterion, in analogy with Proposition \ref{prop:curve_testing}, holds only in the setting of algebraic $\Dscr_X$-modules for a smooth algebraic variety $X$. Then a $\Dscr_X$-module $M$ is regular if and only if $i^*M$ is regular for every locally closed embedding $i:C\to X$ of a smooth algebraic curve $C$.
}
 and holonomic cohomology.  Roughly these conditions guarantee that the solution spaces are finite dimensional and have moderate growth. We denote the derived category of bounded complexes of $\Dscr_X$ modules with regular and holonomic cohomology by $D^b_{rh}(\Dscr_X)$. On the topological side, we look at those bounded complexes of sheaves of $\CC$-vector spaces on $X$ that have constructible cohomology, and we denote the corresponding derived category by $D^b_c(\CC_X)$. 

Under these restrictions the derived solution functor gives the celebrated anti\-equivalence of categories. 

\begin{theorem}\cite{Kas6,Kas10,Me4,BB}\footnote{Beilinson and Bernstein proved the theorem in the algebraic setting.}
There is an anti-equivalence of  (triangulated) categories
\[
R\uHom_{\Dscr_X}(-,\Oscr_X): D_{rh}^{b}(\Dscr_X)^{o}\xrightarrow{\sim}D_{c}^{b}(\CC_X).
\]
\end{theorem}

First we remark that we really need to pass to the derived level contrary to Deligne's Riemann-Hilbert correspondence. Indeed, as mentioned earlier the solution functor is not exact so it cannot induce an equivalence of abelian categories. This theorem is from an algebraic point of view a real advancement, and a vast generalisation of Deligne's Riemann Hilbert correspondence, since we can in particular to every (regular holonomic) $\Dscr_X$-module associate a topological object, a complex of sheaves of $\CC$-vector spaces on $X$ (with constructible cohomology), and vice versa. 

These associated complexes are also rather special, they form an abelian category, and they are called {\em perverse sheaves}, i.e. 
\[
{\rm Perv}(X):=R\uHom_{\Dscr_X}(-,\Oscr_X)({\rm mod}_{rh}\,\Dscr_X^{o})[\dim X].
\]

\section{Homological mirror symmetry symmetries}
We divert the story to mirror symmetry. There we will encounter representations of some fundamental groups and our aim will be to realise them as monodromy representations of differential equations. 

\subsection{Mirror symmetry}
Let us first very briefly say a few words on mirror symmetry, a theory that has its origins in physics, more precisely in string theory.  Typically the spaces that appear in this context have both a complex and a symplectic structure. Moreover, the spaces come in mirror pairs $X$ and $X^o$, with the complex and symplectic structures interlaced. The complex geometry of $X$ mirrors the symplectic geometry of its mirror $X^o$, and vice versa. The picture is still highly speculative. We refer to \cite[Introduction]{MSbook} for a survey of its origins and multiple predictions that mirror symmetry provides to algebraic geometry.

\subsection{HMS categorical symmetries}
Mirror symmetry has been enhanced to a homological statement about the equivalence of certain categories (the derived category and the Fukaya category) that reflect complex and symplectic geometry, respectively. 
The correspondence has been conjectured by Kontsevich \cite{Kontsevich} and nowadays goes under the name of {\em homological mirror symmetry}. 

We discuss here one of the consequences of homological mirror symmetry. For a more precise  explanation of heuristics see \cite[\S 1.1]{HLSam}.
Assume that we regard $X$ as a complex manifold. Then the symplectic structure of the mirror $X^o$ is fixed, but there is still  room for different complex structures. Denote by $\Kscr_X$ the space of complex structures of $X^0$.\footnote{$\Kscr_X$ is also called the ``stringy Kähler moduli space'' (SKMS) of $X$ (i.e. the space of Kähler structures on $X$ coming from symplectic geometry of $X$). 
 The tangent space to the SKMS is $H^2(X,\CC)$ (the space of complexified symplectic forms). However,  there is no global definition, $\Kscr_X$ has only been explicitly defined in very few examples, the difficulty being the determination of the mirror pair.}

Then homological mirror symmetry predicts the following. 

\begin{conjecture}\label{conj}
There exists an action\footnote{We might think of $D^b(X)$ as bounded complexes of vector bundles on $X$.}
\[
\pi_1(\Kscr_X)\curvearrowright D^b( X).
\]
\end{conjecture}

As an immediate corollary of this we would get the following result about the Grothendieck group of $X$. 

\begin{corollary}\label{cor}
There exists an action
\[
\pi_1(\Kscr_X)\curvearrowright K_0( X)_\CC.
\]
\end{corollary}
It is this action about which we will wonder which system of differential equations it corresponds to. 

\subsection{Example}\label{sec:ex}
We look at the conifold,  
$Y={\rm Spec} (\CC[x,y,z,u]/(xu-yz))$.\footnote{One can describe the conifold also as a cone over ${\mathbb P}^1\times {\mathbb P}^1$.} We define
$X={\rm Bl}_{(x,y)}Y$, a small resolution of $Y$. ({In the framework of toric geometry we might represent $Y$ as a cone in $\RR^3$ over the unit square in $\RR^2\times\{1\}$. To obtain $X$ we should add a diagonal hyperplane.})

There is another viewpoint that will be more useful for us. Let $\CC^*$ act on $\CC^4$ as $t.(v_1,v_2,v_3,v_4)=(t^{-1}v_1,t^{-1}v_2,t v_3,t v_4)$. Then we may view 
$Y$ as the (categorical) quotient $\CC^4\quot \CC^*$ ($=\Spec[\CC[x_1,x_2,x_3,x_4]^{\CC^*}$, $x_i=v_i^*$).\footnote{The homomorphism  $\CC[x,y,z,u]/(xu-yz)\to \CC[x_1,x_2,x_3,x_4]^{\CC^*}$, $x\mapsto x_1x_3,y\mapsto x_1x_4,z\mapsto x_2x_3,u\mapsto x_2x_4$, is an isomorphism.} We obtain $X$ as the (GIT) quotient $(\CC^4\setminus V(x_1,x_2))\quot \CC^*$.\footnote{Let as assume that $t.v=tv$, and take $s=v^*$, assume that $\deg s=1$, $\deg x_i=0$, $1\leq i\leq 4$. Then the  GIT quotient $(\CC^4\setminus V(x_1,x_2))\quot \CC^*$ is defined as $\Proj (\CC[x_1,x_2,x_3,x_4,s]^{\CC^*})$.\label{foot:git}}

\medskip

Heuristics from physics \cite{Aspinwall} yield  that $\Kscr_X=\PP^1\setminus \{0,1,\infty\}$. 

To construct a representation of $\pi_1(\Kscr_X)$ on $D^b(X)$ we first view $D^b(X)$ as the (full thick) subcategory of $D^b([\CC^4/\CC^*])$\footnote{Here $[\CC^4/\CC^*]$ denotes the quotient stack. The category $\mod([\CC^4/\CC^*])$ consists of $\CC^*$-equivariant $\CC[x_1,x_2,x_3,x_4]$-modules and the category $\coh([\CC^4/\CC^*])$ of $\CC^*$-equivariant coherent sheaves on $\CC^4$. It follows $D^b([\CC^4/\CC^*])=D^b(\mod([\CC^4/\CC^*]))=D^b(\coh([\CC^4/\CC^*]))$.}, generated by $\Oscr_{\CC^4}$, $\Oscr_{\CC^4}\otimes V(1)$, where $V(n)$ denotes the irreducible (1-dimensional) representation of $\CC^*$ with character $n$, i.e. $t.v=t^nv$ for $v\in V(n)$, \cite[Theorem 8.6]{VdB04}. 

Then it turns out that in the basis $\{\Oscr_{\CC^4}\otimes V(1)$, $\Oscr_{\CC^4}\}$ the action of the three generating loops $\gamma_0,\gamma_1,\gamma_\infty\in \pi_1(\Kscr_X)$ is given by 
\[
\gamma_1=
\left(\begin{matrix}
1&0\\
0&1\\
\end{matrix}\right),
\quad
\gamma_0=
\left(\begin{matrix}
2&1\\
-1&0\\
\end{matrix}\right),
\quad
\gamma_\infty=
\left(\begin{matrix}
0&-1\\
1&2\\
\end{matrix}\right).
\]
See e.g. \cite{HLSam,Donovan,SVdB10}. 

\section{HMS symmetries: toric varieties}
We will approach the conjecture in the setting of toric varieties.

\subsection{Setting}\label{sec:setting}
We assume that $W=\CC^d$ is a $T:=(\CC^*)^n$-representation which is unimodular (i.e. the sum of weights is equal to $0$). 

We describe how to obtain an analogue of the variety $X$ in the case of the conifold, cf. \S\ref{sec:ex}. We should remove some undesirable locus of $W$ and then take the (GIT) quotient. The variety (or stack) $X$ that we obtain in this way is a (crepant) resolution of singularities of $W\quot T$ ($=\Spec\CC[W]^T$). 

Let $X(T)$ be the character group of $T$ and $Y(T)$  the group of $1$-parameter subgroups of $T$. We take a generic $\chi\in X(T)_\RR$.
Let $W^{\chi,u}$ be the $\chi$-unstable locus, i.e.  the set of points $w\in W$ such that if  $\lim_{t\to 0}\lambda(t)w$ for $\lambda\in Y(T)$ exists then $\chi(\lambda)\geq 0$. Then we take 
\[X=[(W\setminus W^{\chi,u})/T].\]
 This is a priori a Deligne-Mumford quotient stack, a quotient stack whose points have finite stabilizers. In the case that all stabilizers are trivial, the corresponding (GIT) quotient variety can replace the stack (i.e. in this case the quotient stack and the quotient variety are isomorphic). The GIT quotient is defined in the analogy with Footnote \ref{foot:git}.\footnote{We assume that $V=\CC v$ is the $1$-dimensional $T$-representation with character $\chi$, i.e. $t.v=\chi(t)v$. Let $w_i$ be a basis of $W$ such that $t.w_i=\beta_i(t)w_i$ for $\beta_i\in X(T)$. Set $x_i=w_i^*$, $1\leq i\leq d$, $d=v^*$. 
We assume that $\deg x_i=0$ and $\deg s=1$. Then  $(W\setminus W^{\chi,u})\quot T:=\Proj(\CC[x_1,\dots,x_d,s]^T)$.}

%It turns out that $X$ has finite quotient singularities. In case that it is not smooth we should replace it with the associated quotient stack.\footnote{We should take $[(W\setminus W^\chi_{ns})/T]$, where $W^\chi_{ns}$ is the $\chi$-unstable locus, i.e.  the set of points $w\in W$ such that if  $\lim_{t\to 0}\lambda(t)w$ for $\lambda\in Y(T)$ exists then $\chi(\lambda)\geq 0$.}

\begin{remark}
The varieties above are exactly affine normal Gorenstein toric varieties whose class group is a torus (i.e. it has no finite group part).
\end{remark}

\subsection{Space of complex structures on $X^o$}
In the case of toric varieties physics heuristics are  rather reliable. In \cite[\S4.1]{DonovanSegal} there is an explicit recipe for $\Kscr_X$ that refers for evidence to \cite{CCIT}.\footnote{The heuristics are  derived from the speculations that a mirror is given by a family of Landau-Ginzburg models \cite{HoriVafa}. See also \cite{Iritani,MirorrSymmetryConstructions}.}

%We introduce a bit of notation that we might need later on. We set $d=\dim W$ and denote $\TT=(\CC^*)^d$. Recall that $T=(\CC^*)^n$ and that $(\beta_i)_{i=1}^d$ are ($T$-)weights of $W$.  

Set $d=\dim W$. Let $(\beta_i)_{i=1}^d$ be $T$-characters of $W$. Note that $X(T)\cong \ZZ^n$ and set $B=(\beta_i)_{i=1}^d\in M_{n\times d}(\ZZ)$.  We define $A$ (up to an automorphism of $\ZZ^{d-n}$) by the following exact sequence
\begin{equation}\label{eq:exseq}
0\to \ZZ^{d-n}\xrightarrow{A}\ZZ^d\xrightarrow{B} \ZZ^n\to 0.
\end{equation}
Then $\Kscr_X$ is the complement of a hypersurface $V(E_A)\subset T$ where $E_A$ is the {\em principal $A$-determinant}. We refer  to \cite[\S10.1.A]{GKZbook} for the definition.\footnote{In loc.cit. $E_A$ stands for $A'$ where $A=(A',1)$ which we may assume since $\sum_i\beta_i=0$.} Alternatively, see \cite{DonovanSegal,Kite}.

In a sufficiently symmetric case $V(E_A)$ is much simpler.

\begin{theorem}\cite{Kite}\label{thm:Kite}
If $W$ is quasi-symmetric\footnote{$W$ is quasi-symmetric if for all lines $0\in\ell\in X(T)_\RR$, $\sum_{\beta_i\in\ell}\beta_i=0$.} then $\Kscr_X$ is the complement of a hyperplane arrangement (in logarithmic coordinates) in $T=(\CC^*)^n$.
\end{theorem}

The hyperplane arrangement in $1/(2\pi i)\log T=X(T)_\CC$ can be explicitly described. Let $\Delta$ be the Minkow\-ski sum of $[0,(1/2)\beta_i]$. Let $(H_i)_i$ be the supporting (affine) hyperplanes of $\Delta$. Then the hyperplane arrangement is the complexification of the real hyperplane arrangement $\cup_i (-H_i)+X(T)$ (up to a suitable translation). This is an infinite, but locally finite hyperplane arrangement.

This hyperplane arrangement was prior to the result of Kite heuristically predicted to coincide with $\Kscr_X$ in \cite{HLSam}.

\begin{example}
We make a quick sanity check in the case of the conifold, cf. \S\ref{sec:ex}. Then $\Kscr_X=\PP^1\setminus \{0,1,\infty\}$. Applying $1/(2\pi i)\log$ to $\PP^1\setminus\{0,1,\infty\}=\CC\setminus\{0,1\}$ we obtain $\CC\setminus \ZZ$.  On the other hand, by the above recipe, $\Delta=[-1,1]$ (as $(\beta_i)_{i=1}^4=(-1,-1,1,1)$) and the hyperplane arrangement is given by $\ZZ$. Thus, the two descriptions are consistent. 
\end{example}

\subsection{HMS symmetries: quasi-symmetric case}
Assume that $\CC^d$ is a quasi-symmetric representation of $(\CC^*)^n$. In this case Halpern-Leistner and Sam \cite{HLSam} confirmed Conjecture \ref{conj}. 

\begin{theorem}\cite{HLSam}\label{thm:HLS}
There exists an action of $\pi_1(\Kscr_X)$ on $D^b(X)$.
\end{theorem}

As in \S\ref{sec:ex}, $D^b(X)$ is identified with the (full thick) subcategory $D$ of $D^b([W/T])$ generated by $\{\Oscr_W\otimes V(\mu)\mid \mu\in (\nu+\Delta)\cap X(T)\}$, where $V(\mu)$ is the irreducible $T$-representation with character $\mu$, and $\nu\in X(T)_\RR$ is generic \cite{HL,SVdB}\footnote{$\nu$ is not parallel to any face of $\Delta$.}.

Then this action can be explicitly described, especially relying on the concrete description of the fundamental group of the complement of a complexified hyperplane arrangement \cite{Salvetti}.
% (and the combinatorial algorithm \cite{SVdB} for relating different $\Oscr_W\otimes V(\mu)$.). 
See \S\ref{subsubsec:psch} below.

\begin{remark}
The statement can be generalised to some reductive groups, i.e. those groups $G$ for which $X(G)
%=X(T)^\Wscr
\neq 0$,  
%where $\Wscr$ denotes the Weyl group of $G$, 
if some genericity assumptions are satisfied.\footnote{The condition $\sum_i \RR\beta_i=X(T)$ should be satisfied and there  should exist $\chi\in X(G)$ which is not parallel{} to any face of $\Delta$.} See \cite{HLSam}.
\end{remark}

\section{HMS differential equations: quasi-symmetric case}
In this section we assume that we are in the setting of \S\ref{sec:setting}. Moreover, we assume that $W$ is quasi-symmetric. Having Theorem \ref{thm:HLS}, providing evidence for Conjecture \ref{conj}, at our disposal, we also obtain Corollary \ref{cor}. Hence, $\pi_1(\Kscr_X)$ acts on $K_0(X)_\CC$. We want to determine which (regular) system of differential equations on $(\CC^*)^n$ this action corresponds to.% cf. Corollary \ref{cor:Deligne}. 

\subsection{Example}\label{subsec:ex}
We first want to understand the monodromy representation in the case of the conifold, cf. \S\ref{sec:ex}.

We look at the {\em Gauss hypergeometric equation}
\[
z(1-z)y''+(c-(a+b+1)z)y'-ab y=0.
\]
The monodromy is given by, see e.g. \cite{BeukersHeckman},
\begin{align*}
&\begin{aligned}
\gamma_1&=
\left(\begin{matrix}
1&-e^{2\pi \mathrm{i}(c-b)}-e^{2\pi \mathrm{i}(c-a)}+e^{2\pi \mathrm{i} c}+1\\
0&e^{2\pi \mathrm{i}(c-a-b)}\\
\end{matrix}\right),
\\
\gamma_0&=
\left(\begin{matrix}
1+e^{-2\pi \mathrm{i} c}&1\\
-e^{-2\pi \mathrm{i} c}&0\\
\end{matrix}\right),\\
\gamma_\infty&=
\left(\begin{matrix}
0&-e^{2\pi \mathrm{i} (a+b)}\\
1&e^{2\pi \mathrm{i} a}+e^{2\pi \mathrm{i} b}\\
\end{matrix}\right).
\end{aligned}&
\end{align*}

Setting $a=b=c=0$ we obtain matrices that we have already encountered in \S\ref{sec:ex}. From this one may deduce that the action of $\pi_1(\Kscr_X)$ on $K_0(X)_\CC$ from Theorem \ref{thm:HLS} in the case of the conifold corresponds to $z(1-1)y''-zy'=0$, i.e. the Gauss differential equations with parameters $a=b=c=0$ (which is a regular on $\PP^1$ with singularities at $0,1,\infty$).

\subsection{Example with parameters}\label{subsec:expar}
We change the focus a bit, and ask whether we can find an action of $\pi_1(\Kscr_X)$ on $K_0(X)_\CC$ that would  give the Gauss hypergeometric equation also for other parameters. We obtained the original action from an action of $\pi_1(\Kscr_X)$ on $D^b(X)$. We would want to tweak this action a little bit to open the route to other parameters. 

For this first observe that $(\CC^*)^4$ acts on $\CC^4$ coordinate-wise. The initial $\CC^*$ embeds in it via the map $t\mapsto (t^{-1},t^{-1},t,t)$ determined by the action of $\CC^*$ on $\CC^4$, cf. \S\ref{sec:ex} and \eqref{eq:exseq}. This inclusion splits, and the complement is $(\CC^*)^3$. We seem to be well on the way, the dimension of the complement torus coincides with the number of  parameters in the Gauss hypergeometric equation. 

Now follows a slightly more technical part. To get an action for other $a,b,c$ we need to replace $D^b(X)$ by a bigger category $\tilde{D}$ such that $X((\CC^*)^3)$ acts on it. 

We define $\tilde{D}$ as the (full thick) subcategory of $D^b([\CC^4/(\CC^*)^4]$ generated by $\Oscr_{\CC^4}\otimes V(\mu)$, $\mu\in X((\CC^*)^4)$ such that $B\mu\in \{0,1\}$ (see \eqref{eq:exseq} for $B$). 

It turns out that $\pi_1(\Kscr_X)$ still acts on $\tilde{D}$. However, $K_0(\tilde{D})_\CC$ is a (free rank $2$) module over $\CC\{X((\CC^*)^3)\}\cong \CC[(\CC^*)^3]$\footnote{Here $\CC\{X((\CC^*)^3)\}$ is the group algebra of $X((\CC^*)^3)$, while $\CC[(\CC^*)^3]$ is the coordinate ring of $(\CC^*)^3$.}. Specialising at (sufficiently generic\footnote{This is in particular satisfied if $a,b,a-c,b-c$ are all non-integers. However, one might check that $a=b=c=0$ as in \S\ref{subsec:ex} also work.}) $h\in (\CC^*)^3$ we obtain an action of $\pi_1(\Kscr_X)$ on a $2$-dimensional $\CC$-vector space. This action corresponds to the Gauss hypergeometric equation with parameters $-1/(2\pi i) \log h$.

\subsection{GKZ hypergeometric systems}
The GKZ hypergeometric systems are systems of differential equations that generalise the Gauss hypergeometric differential equation, as well as Appell, Lauricella, Horn... They were introduced and studied by Gelfand, Kapranov and Zelevinsky \cite{GFK0,GKZhyper,GKZEuler}. Allegedly, they were introduced as a unified approach to the multidimensional generalisations of the Gauss hypergeometric functions. In some sense, the construction of the GKZ hypergeometric system is dictated by the desired set of solutions, which should be hypergeometric power series. See Remark \ref{rem:pssol}.
%The suitable generalisation of the Gaussian hypergeometric equation suitable for other torus representations is given 

Let $\alpha\in \CC^{d-n}$. Recall the exact sequence \eqref{eq:exseq}. Let $B^*:\ZZ^n\to \ZZ^d$ be the dual of $B$. 
%We may write it more canonically, using $T$, $\TT$. Additionally, we define $H=\coker(T\hookrightarrow \TT)$. 
Then the hypergeometric GKZ system with parameter $\alpha$ is defined by the following differential operators:
\begin{itemize}
	\item homogeneity relations: $\sum_{j=1}^d a_{ij}x_j\partial_j-\alpha_i$, $1\leq i\leq d-n$,
	\item box relations: $\square_l=\prod_{l_i>0}\partial_i^{l_i}-\prod_{l_i<0}\partial_i^{-l_i}$, $l\in B^*\ZZ^n$.
\end{itemize}

Note that this is a system of differential equations on $(\CC^*)^d$. However, the homogeneity relations allow to descend these differential equations to $(\CC^*)^n$.\footnote{The corresponding $D$-module on $(\CC^*)^d$ is weakly equivariant for the action of $(\CC^*)^{d-n}$, hence it descends to $(\CC^*)^n$, see e.g. \cite[Corollary A.11]{SVdBgkz}.}

This descent also allows us to recover the Gauss hypergeometric equation from the GKZ hypergeometric system corresponding to the conifold, i.e. for $B=(-1,-1,1,1)$. 
\begin{example}
In the case of the conifold, we may take 
\[
A=\left(
\begin{matrix}
-1&1&0&0\\
1&0&1&0\\
1&0&0&1\\
\end{matrix}
\right).
\]
Then  a solution $\Phi$ of the GKZ hypergeometric system satisfies
\begin{align*}
(-x_1\partial_1+x_2\partial_2)\Phi&=\alpha_1\Phi\\
(x_1\partial_1+x_3\partial_3)\Phi&=\alpha_2\Phi\\
(x_1\partial_1+x_4\partial_4)\Phi&=\alpha_3\Phi\\
\partial_1\partial_2-\partial_3\partial_4&=0.
\end{align*}
Setting $\alpha=(c-1,-a,-b)$ a simple manipulation yields
\[
(x_3^{-1}x_4^{-1}(x_1\partial_1^2-(1+a+b)x_1\partial_1-ab)-x_2^{-1}(x_1\partial_1^2-c\partial_1))\Phi=0.
\]
Then $F(x):=\Phi(x,1,1,1)$ is a solution of the Gauss hypergeometric equation. Moreover, by homogeneity relations, $F$ determines $\Phi$. 
\end{example}

We denote the corresponding $\Dscr_{(\CC^*)^d}$-module, cf. \S\ref{sec:Dmodules}, by $\Pscr(\alpha)$, and its restriction to $(\CC^*)^n$  by $P(\alpha)$. The next proposition reveals that they are well-behaved, as required for the Riemann-Hilbert correspondence.

\begin{proposition}\cite{Adolphson}
The $\Dscr_{(\CC^*)^d}$-module $\Pscr(\alpha)$ is holonomic with regular singularities. The same holds for the $\Dscr_{(\CC^*)^n}$-module $P(\alpha)$.
\end{proposition}

\begin{remark}\label{rem:pssol}
To follow on the introduction to this subsection we record here that the multidimensional hypergeometric (formal) series\footnote{We abuse the notation and denote by $B^*$ also the ``complexified'' $B^*:\CC^n\to \CC^d$.}
\[
\Phi_\gamma(x_1,\dots,x_d)=\sum_{l\in \ZZ^n}\prod_{i=1}^d \frac{x_i^{B^*l+\gamma_i}}{\Gamma(B^*l+\gamma_i+1)},
\]
where $\gamma\in (\CC)^d$ is such that $A\gamma=\alpha$, is a (formal) solution of the GKZ hypergeometric system \cite{GKZhyper}.\footnote{By appropriately varying $\gamma$ one can achieve that such power series a basis of solutions that converge on an open set.} 
Moreover, also Euler integrals generalise \cite{GKZEuler} to give solutions to the GKZ hypergeometric system. Another handy class of solutions is given by Mellin-Barnes integrals \cite{Beukers} that we crucially employ in the proof of Theorem \ref{thm:SVdB} below. 
\end{remark}

\subsection{Decategorification of HMS symmetries}
We want to determine the system of differential equations whose monodromy representation coincides with the representation of $\pi_1(\Kscr_X)$ on $K_0(X)_\CC$ obtained from Theorem \ref{thm:HLS}. However, we cannot quite do that. Instead, we tweak the action a bit, as in \S\ref{subsec:expar}.

Analogously to \S\ref{subsec:expar} we note that $(\CC^*)^d$ acts on $\CC^d$ coordinate-wise, and we have the inclusion $T=(\CC^*)^n\hookrightarrow (\CC^*)^d$ which splits. The complement is $(\CC^*)^{d-n}$.  

We replace $D^b(X)$ by a bigger category $\tilde{D}$, the (full thick) subcategory of $D^b([\CC^d/(\CC^*)^d])$ generated by 
\[
\left\{\Oscr_{\CC^d}\otimes V(\mu)\mid B\mu\in (\nu+\Delta)\cap X(T)\right\},
\]
 cf. the paragraph below Theorem \ref{thm:HLS}. Then $X((\CC^*)^{d-n})$ acts on $\tilde{D}$ and $K_0(\tilde{D})_\CC$ is a $\CC\{X((\CC^*)^{d-n})\}\cong\CC[(\CC^*)^{d-n}]$-module. 

\begin{theorem}\cite{SVdBgkz}\label{thm:SVdB}
Assume that $\alpha\in \CC^{d-n}$ is generic.\footnote{We require that $\alpha$ is non-resonant, i.e. $\alpha$ does not belong to the hyperplane arrangement consisting of $\ZZ^{d-n}$-translates of the supporting hyperplanes of the cone $\RR_+A$.}  
The monodromy representation of the GKZ system of differential equations with parameter $\alpha$ restricted to $\Kscr_X$ is isomorphic to the representation of $\pi_1(\Kscr_X)$ on $K_0(\tilde{\Dscr})_\CC$ specialised at $e^{-2\pi i\alpha}$\footnote{More precisely, $K_0(\tilde{\Dscr})_\CC\otimes_{\CC[(\CC^*)^{d-n}]} \CC$ for $\CC[(\CC^*)^{d-n}]\to \CC$, $p\mapsto p(e^{-2\pi i\alpha})$.}.
\end{theorem}

As a corollary we  obtain in particular a description of the full monodromy of such ``quasi-symmetric'' GKZ hypergeometric systems. In \cite{Beukers}, Beukers  describes the ``local'' monodromy.

\begin{remark}
There are various other results where an interesting system of differential equations is obtained from  actions on derived categories, often also inspired by mirror symmetry. We mention here \cite{BHMB,ABM,BridgelandYuSutherland}.
\end{remark} 

\section{Liftings}
Theorem \ref{thm:HLS} (and accordingly Theorem \ref{thm:SVdB}) extend a bit further, in analogy with D-modules introduced in \S\ref{sec:Dmodules} and the associated perverse sheaves, defined as the image of the abelian category of D-modules by the derived solution functor \S\ref{subsec:RH}. 

\subsection{Perverse schobers}
Recall that a representation of $\pi_1(\Kscr_X)$ corresponds to a local system on $\Kscr_X$, cf. \ref{subsec:DRH}. If $\pi_1(\Kscr_X)$ acts instead   on a category, we might say that is corresponds to  a local system of categories on $\Kscr_X$. In the quasi-symmetric setting, $\Kscr_X$ is  a complement of a hyperplane arrangement in $(\CC^*)^n$ (in logarithmic coordinates), cf. Theorem \ref{thm:Kite}. We may extend a local system on $\Kscr_X$ to a perverse sheaf on $(\CC^*)^n$. This extension for the particular action of Theorem \ref{thm:HLS} also lifts on the level of derived categories, and we get what we might call a perverse sheaf of categories on $(\CC^*)^n$ \cite{SVdB10}. It also goes under the name of a {\em perverse schober}, which was coined by Kapranov and Schechtman \cite{KapranovSchechtmanSchobers} for a categorification of a perverse sheaf.

The rest of this subsection builds on this extension, and in return also illuminates the proof of Theorem \ref{thm:HLS}. Unfortunately, it is rather technical.

\subsubsection{Perverse sheaves over real hyperplane arrangements}\label{subsubsec:ps}
While in general the abelian category of perverse sheaves might be difficult to describe, in the case of complements of complexified real hyperplane arrangements there exists a concrete combinatorial description \cite{KapranovSchechtman}, which is apt for categorification. 

Let $\Hscr$ be an affine hyperplane arrangement in a finite dimensional real vector space $V=\RR^n$.
%The closures of the connected components of $\RR^n\setminus \Hscr$ are convex polytopes. 
Then $\Hscr$ stratifies $V$ into a set $\Cscr$ of locally closed subsets.\footnote{The elements of $\Cscr$ are level sets for  $({\rm sign}\, f_H)_{H\in \Hscr}$ where $f_H$ is the affine map defining $H$.}
We partially order $\Cscr$ by 
  $C'\le C$ iff $C'\subset \overline{C}$. 
A triple of faces $(C_1,C_2,C_3)$ is {\em collinear} if there exists $C'\leq C_1,C_2,C_3$ and there exist $c_i\in C_i$ such that $c_2\in[c_1,c_3]$.

We denote by ${\rm vec(\CC)}$ the category of finite dimensional $\CC$-vector spaces.

\begin{theorem}\cite{KapranovSchechtman}\label{thm:psrh}
  The category %$\Perv(V_\CC,\Hscr_\CC)$ 
   of perverse sheaves  on $V_\CC$ with respect to the stratification induced by $\Hscr_\CC$ is equivalent to the category
  of diagrams consisting of
\begin{itemize}
	\item
   finite dimensional vector spaces $E_C$,
  $C\in\Cscr$, and 
  \item
  linear maps $\gamma_{C'C}:E_{C'}\to E_C$,
  $\delta_{CC'}:E_C\to E_{C'}$ for $C'\leq C$
  \end{itemize} 
  such that
  $(E_C,(\gamma_{C'C})_{CC'})$, resp. $(E_C,(\delta_{CC'})_{CC'})$, is a
  representation of $(\Cscr,\leq)$, resp. $(\Cscr,\geq)$, in ${\rm vec}(\CC)$, and the following
  conditions are satisfied:
\begin{itemize}
\item $\gamma_{C'C}\delta_{CC'}=\id_{E_C}$ for $C'\leq C$. 
In particular, $\phi_{C_1C_2}:=\gamma_{C'C_2}\delta_{C_1C'}$ for $C'\leq C_1,C_2$ is well defined.
\item $\phi_{C_1C_2}$ is an isomorphism for every $C_1\neq C_2$ of the same dimension $\ell$, which lie in the same $\ell$-dimensional
  affine space and share a facet. 
\item $\phi_{C_1C_3}=\phi_{C_2C_3}\phi_{C_1C_2}$ for collinear
  triples of faces $(C_1,C_2,C_3)$.
\end{itemize}
\end{theorem}

\subsubsection{Perverse schobers over real hyperplane arrangements}\label{subsubsec:psch}
%\subsubsection{}
To define perverse schobers over real hyperplane arrangements we may word for word translate the description of perverse sheaves from Theorem \ref{thm:psrh} to the setting of triangulated categories. When we apply $K_0(-)_\CC$ we  get back the data defining a perverse sheaf.

\begin{definition}\cite{BondalKapranovSchechtman}\label{def:schober}
  A perverse schober  on $V_\CC$  with respect to the stratification induced by $\Hscr_\CC$\footnote{A perverse schober in this context is also called an $\Hscr$-schober.} is given by 
\begin{itemize}
	\item
  triangulated categories $\Eescr_C$, 
  $C\in \Cscr$,  and
  \item
  adjoint exact functors $(\delta_{CC'}:\Eescr_C\to \Eescr_{C'},\gamma_{C'C}:\Eescr_{C'}\to \Eescr_C)$ 
   for $C'\leq C$
\end{itemize}
    such that
$(\Eescr_C,(\delta_{C'C})_{C'C})$ defines a pseudo-functor from $(\Cscr,\ge )$ to the 
 $2$-category of
  triangulated categories satisfying  the following conditions:
 \begin{itemize}
 \item The unit of the adjunction $(\delta_{CC'},\gamma_{C'C})$
defines a natural isomorphism
   $\id_{\Eescr_C}\xrightarrow{\cong}\gamma_{C'C}\delta_{CC'}$ for $C'\leq C$ , and
   thus $\phi_{C_1C_2}:=\gamma_{C'C_2}\delta_{C_1C'}$ for
   $C'\leq C_1,C_2$ is well defined up to canonical
natural   isomorphism.
\item
$\phi_{C_1C_2}$ is an equivalence for every $C_1\neq C_2$ of the same dimension $\ell$, which lie in the same $\ell$-dimensional
  affine space and share a facet. 
\item
 The counit of the adjunction $(\delta_{C_0C_2},\gamma_{C_2C_0})$
 defines a natural isomorphism
$\phi_{C_2C_3}\phi_{C_1C_2}\xrightarrow{\cong} \phi_{C_1C_3}$ for collinear triples of faces $(C_1,C_2,C_3)$.
\end{itemize}
\end{definition}

This definition also sheds some light on the proof of Theorem \ref{thm:HLS}, cf. the paragraph following it, and allows its extension. 

\begin{theorem}\cite{SVdB10}
The local system on $\Kscr_X$ from Theorem \ref{thm:HLS} extends to a perverse schober on $(\CC^*)^n$.\footnote{We identify  $(\CC^*)^n$ with $\CC^n/\ZZ^n$, and in order to use Definition \ref{def:schober} we should also impose an action of $\ZZ^n$ on a perverse schober, which consists of isomorphisms $\Escr_C\to \Escr_{gC}$ for $g\in \ZZ^n$ satisfying some compatibility conditions, see e.g. \cite[\S3.3]{SVdB10}.}
\end{theorem}

\begin{remark}
By a suitable tweak as in Theorem \ref{thm:SVdB}, we obtain perverse schobers whose decategorifications are the perverse sheaves obtained as solution complexes of GKZ hypergeometric D-modules \cite{SVdBgkz}. 
\end{remark}

\subsection{HMS predictions}
GKZ hypergeometric systems appeared here rather adhoc, and not really motivated. In fact, it is homological mirror symmetry that indicates that they should be there \cite{Iritani,CCITHodge,BHMB}. 

While we only combinatorially match the two perverse sheaves, one would desire to construct a canonical correspondence via the following sequence of maps (GM denotes Gauss-Manin)
\begin{multline*}
(K_0(D^b([W/T]))\supset)\, K_0(D) \xrightarrow{\sim} K_0(X)\xrightarrow{\sim} \\H^*(X) \;
\text{(for. quantum conn.)}
\xrightarrow{\text{mirror map}}
\text{\{rel. tw. DR-coh. at $\infty$\} (for. GM conn.)}\\\xrightarrow{\text{anal. cont.}} \{\text{solutions to GKZ system}\}.
\end{multline*}

However,  the heuristics of why this action would lift to an action on the derived category of $X$ are still somewhat mysterious.\footnote{The RHS of the mirror map corresponds to the $B$-side of the LG-model, which would in turn lead to a GKZ system on the Fukaya category of $X$, rather than on the derived category.}

\section{Acknowledgement}
We are foremost grateful to Michel Van den Bergh for a journey to kaleidoscopic areas that would otherwise remain inaccessible to us. 

Moreover, we thank Geoffrey Janssens, Urban Jezernik, Igor Klep and the referee for a generous assortment of comments and suggestions. 
%\bibliography{nccr,ecm}
\bibliographystyle{amsalpha} 

\end{document}